\documentclass[reqno, 12pt]{amsart}
\usepackage{amsmath}
\usepackage{amssymb}
\usepackage{amsthm}
\usepackage{mathtools}
\usepackage{tikz-cd}
\usepackage{enumerate}
\usepackage[mathscr]{eucal}
\usepackage{graphicx}
\usepackage{epsfig}
\newtheorem{theorem}{Theorem}[section]

\newtheorem{proposition}[theorem]{Proposition}
\newtheorem{corollary}[theorem]{Corollary}
\newtheorem{definition}[theorem]{Definition}
\theoremstyle{definition}
\newtheorem{example}[theorem]{Example}
\newtheorem{remark}[theorem]{Remark}

\begin{document}
	\title[]{Hyperspace convergences, bornologies and geometric set functionals}
	\author{Yogesh Agarwal \and Varun Jindal}
	
	\address{Yogesh Agarwal: Department of Mathematics, Malaviya National Institute of Technology Jaipur, Jaipur-302017, Rajasthan, India}
	\email{yagarwalm247@gmail.com}
	
	\address{Varun Jindal: Department of Mathematics, Malaviya National Institute of Technology Jaipur, Jaipur-302017, Rajasthan, India}
	\email{vjindal.maths@mnit.ac.in}
	
	\subjclass[2010]{Primary ; Secondary }	
	\maketitle
	\begin{abstract}
		For a bornology $\mathcal{S}$ of subsets of a metric space $(X,d)$, we consider the following unified approaches of hyperspace convergence: convergence induced through uniform convergence of distance functionals ($\tau_{\mathcal{S},d}$-convergence); bornological convergence, and the weak convergence induced by a family of gap and excess functionals.  An interesting problem regarding these convergences is to investigate when any two of them are equivalent.  In this article, we investigate the relation of $\tau_{\mathcal{S},d}$-convergence with the other two convergences, which is not completely transparent. As a main tool for our investigation, we use the idea of pointwise enlargement of a set by a positive Lipschitz function. As applications of our results, we provide new proofs of some known results about Attouch-Wets convergence. 

	\end{abstract}
	
	\section{Introduction}
	A convergence of a net or a sequence of subsets of a topological space is called a hyperspace convergence. One of the widely accepted convergence of sets is the Attouch-Wets convergence (AW-convergence for short) or bounded Hausdorff convergence.  It has proven to be significant for the researchers working in the areas such as optimization, approximation theory,  and variational analysis \cite{attouch1991topology,AW-2,attouch1991quantitative, attouch1993, ToCCoS, beer1993weak, lucchetti2006convexity,penot1991aw, penot2005aw, rockafellar2009variational}. This notion of convergence is specially helpful while dealing with closed convex sets and lower semi-continuous convex functions (when identified with their epigraphs) in the context of normed spaces. 
	Some notable facts in this direction are: convergence of proper lower semi-continuous convex functions in this sense enforces the AW-convergence of their sublevel sets at fixed height greater than the infimal value of the limit (Proposition 7.1.7, \cite{ToCCoS}), also the Fenchel transform is continuous with respect to this convergence \cite{beer1990conjugate}. In particular, one can deduce that the polar map for closed convex sets is continuous with respect to this convergence which suggests the stability of the Attouch-Wets convergence with respect to duality.
	
	The main goal of this paper is to establish connection between various modern (unified) approaches of set convergences studied in the literature. To provide the motivation, we first review various representations of Attouch-Wets convergence. 
	
	For a metric space $(X,d)$, the \textit{Attouch-Wets convergence}  of a net $(A_\lambda)$ to $A$ amounts to the uniform convergence of the net $(d(\cdot, A_\lambda))$ of distance functionals to $d(\cdot, A)$ on bounded subsets of $(X,d)$. This kind of approach is helpful to visualize hyperspace convergence as a convergence in a function space. Another way to look at the AW-convergence is through  geometric set functionals. Two important set functionals are the gap and excess functionals. For a nonempty subset $S$ of a metric space $(X,d)$, the \textit{gap functional} $D_d(S, \cdot)$ on $CL(X)$ (the collection of all nonempty closed subsets of $(X,d)$) is defined as $A \to D_d(S,A)$ for $A \in CL(X)$, where $D_d(S,A) = \inf\{d(x,a): x\in S \text{ and } a \in A\}$. Similarly, the \textit{excess functional} determined by $S$ is given by $A \rightarrow e_d(S, A)= \sup \{d(x,A): x \in S\} \text{ for }A \in CL(X)$. Then the AW-convergence is the convergence corresponding to the weak topology determined by  the family $\{D_d(S,\cdot): S \subseteq X \text{ is bounded}\}\cup \{e_d(S,\cdot): S \subseteq X \text{ is bounded}\}$. In terms of excess functionals, the AW-convergence of a net $(A_\lambda)$ to $A$ means that $\lim e_d(A\cap S, A_\lambda) = 0 = \lim e_d(A_\lambda \cap S, A)$ for each bounded subset $S$ of $X$. These kinds of representations have been significant in epigraphical analysis and regularization of convex functions \cite{ToCCoS}. In the same vein, a net $(A_{\lambda})$ of nonempty closed subsets is AW-convergent to $A$ provided the truncation of $A$ (resp. $A_{\lambda}$) by each nonempty bounded subset of $(X,d)$ is eventually contained in each small enlargement of $A_{\lambda}$ (resp. $A$). This approach has been quite helpful in various scenarios, for example in showing the equivalence of strong convergence of linear functionals and AW-convergence of level sets at fixed heights in the context of normed space \cite{awbornological1}.  
	
	
	If in the above approaches, we replace the family of metrically bounded sets by an arbitrary family $\mathcal{S}$ of nonempty subsets of $X$, then the result is a unified framework to study various hyperspace convergences. Unlike AW-convergence, for an arbitrary family $\mathcal{S}$ the resulting convergences may not be equivalent (\cite{beer2023bornologies,gapexcess, Idealtopo}).  We mainly consider the following three kinds of hyperspace convergences in this paper.
	\begin{enumerate}
		\item The hyperspace convergence through uniform convergence of distance functionals on the members of a family $\mathcal{S}$ (denoted by $\tau_{\mathcal{S},d}$-convergence). This convergence was first formally studied by Beer et al. in \cite{Idealtopo} and also studied in \cite{agarwal2024set,agarwal2024set2,beer2023bornologies}.
		\item The convergence deduced from the weak topology generated by the family $\{D_d(S,\cdot): S \in \mathcal{S}\} \cup  \{e_d(S,\cdot): S \in \mathcal{S}\}$ of gap and excess functionals known as the \textit{Gap and Excess convergence}. The gap and excess convergence for a general family $\mathcal{S}$ was studied by Beer and Levi in \cite{gapexcess}. In the special cases such a convergence also appeared earlier in the literature (see \cite{beer1993weak}).  
		\item The $\mathcal{S}$-convergence of a net $(A_\lambda)$ to $A$, informally, the set convergence obtained by truncating the members of $(A_{\lambda})$ and $A$ by members of $\mathcal{S}$. The $\mathcal{S}$-convergence is also known as bornological convergence was first formally studied by Lechicki et al. in \cite{BCL}.  It is particularly interesting to know when this convergence is topological \cite{beer2009operator, beer2023bornologies,Bcas,gapexcess,PbciAWc, Ucucas, Idealtopo}.         
	\end{enumerate}
	Other than these, the weak topology generated by the family $\{D_d(S,\cdot): S \in \mathcal{S}\}$ of gap functionals is called the \textit{Gap topology} (denoted by $\mathsf{G}_{\mathcal{S},d}$, \cite{beer2013gap}), and by the family  $\{e_d(S,\cdot): S \in \mathcal{S}\}$ of excess functionals is called the \textit{left excess topology} (denoted by $\mathcal{LE}_{\mathcal{S},d}$). The excess topologies for a general family $\mathcal{S}$ have been studied in \cite{beer2014excess}.  Many classical hyperspace topologies can be viewed as gap topologies (see, Corollaries $4.1.6$ and $4.1.7$ in \cite{ToCCoS}) or gap and excess topologies (see, Theorems $4.2.3$ and $4.2.4$ in \cite{ToCCoS}).  	   
	
	
	The central to our study is the $\tau_{\mathcal{S},d}$-convergence for an arbitrary bornology $\mathcal{S}$ on a metric space $(X,d)$ (\cite{ beer2023bornologies, Idealtopo}). Like the other hyperspace convergences, the $\tau_{\mathcal{S},d}$-convergence may be split into two halves: upper ($\tau_{\mathcal{S},d}^{+}$-convergence) and lower ($\tau_{\mathcal{S},d}^{-}$-convergence).   It turns out that for a general bornology $\mathcal{S}$, $\tau_{\mathcal{S},d}$-convergence does not have nice structural properties like the AW-convergence. For instance,  

	\begin{enumerate}
		\item unlike AW-convergence, $\tau_{\mathcal{S},d}$-convergence may not correspond to a metric convergence on $CL(X)$ in general (for example, when $\mathcal{S} = \mathcal{F}(X)$ $\tau_{\mathcal{S},d}$ = Wijsman Convergence);
		\item the $\tau_{\mathcal{S},d}$-convergence does not coincide with bornological convergence in general (\cite{Idealtopo});
		\item  the $\tau_{\mathcal{S},d}$-convergence may not correspond to the weak convergence determined by  gap and excess functionals (Example \ref{counterexampleTsdgapexcess}).
	\end{enumerate}
Thus, the bornology $\mathcal{S}$ needs some additional structure in order for the $\tau_{\mathcal{S},d}$ convergence to have a certain feature.  For example, in \cite{agarwal2024set2}, the authors provide a necessary and sufficient condition on $\mathcal{S}$ for the metrizability of the space $(CL(X),\tau_{\mathcal{S},d})$. Similarly, the second and third problems were partly addressed in \cite{agarwal2024set}, in which the authors gave nice covering and stability conditions on $\mathcal{S}$ for the coincidence of $\tau_{\mathcal{S},d}^{+}$-convergence with the upper gap and upper $\mathcal{S}$-convergences. In \cite{Idealtopo}, the authors investigated the equivalence of upper and lower $\tau_{\mathcal{S},d}$-convergences, respectively with upper and lower $\mathcal{S}$-convergences in terms of convergence of certain nets. 

In the present paper, we will address the following problems.

\begin{itemize}
	\item[(i)] To find the structural property of $\mathcal{S}$ for the coincidence of (lower) $\tau_{\mathcal{S},d}$-convergence and (lower) bornological convergence.
	\item[(ii)] Under what conditions on $\mathcal{S}$, the $\tau_{\mathcal{S},d}$-convergence is the weak convergence determined by gap and excess functionals.\end{itemize}

In particular, we give nice stability and covering conditions on the bornology $\mathcal{S}$ for the equivalence of $\tau_{\mathcal{S},d}$-convergence with the bornological convergence and the gap and excess convergence. Under certain assumption, it is proved that $\tau_{\mathcal{S},d} = \mathcal{S}$ provided $\mathcal{S}$ is stable under Lipschitz enlargements. Similarly, under the same assumption, we show that $\tau_{\mathcal{S},d}$-convergence can be deduced as a weak convergence determined by gap and excess functionals if the Lipschitz enlargements of members of $\mathcal{S}$ fulfill a specific covering condition. New proofs to some known results on Attouch-Wets convergence have been given as corollaries.
	
	
	\section{preliminaries}
	In this section, we define various set convergences and introduce other key concepts. For a metric space $(X,d)$, a nonempty family $\mathcal{S}$ of nonempty subsets of $X$ which is hereditary, stable under finite union, and that forms a cover of $X$ is called a \textit{bornology} on $X$. The smallest bornology on a metric space $(X,d)$ is $\mathcal{F}(X)$, the family of all nonempty finite subsets of $X$, and the largest is $\mathcal{P}_0(X)$, the family of all nonempty subsets of $X$. Some other classical bornologies on a metric space $(X,d)$ are
	\begin{itemize}
		\item $\mathcal{B}_d(X)$ = the family of all nonempty $d$-bounded subsets of $(X,d)$; 
		\item $\mathcal{K}(X)$ = the set of all nonempty relatively compact subsets of $(X,d)$; 
		\item $\mathcal{T}\mathcal{B}_{d}(X)$  = the set of all nonempty $d$-totally bounded subsets of $(X,d)$.
	\end{itemize}  
	
	The open ball centered at $x\in X$ and radius $\epsilon > 0$ in $(X,d)$ is denoted by $B_d(x, \epsilon)$. For $A \subseteq X$ and $\epsilon > 0$, the \textit{$\epsilon$-enlargement of $A$}, denoted by $B_d(A, \epsilon)$, is defined as $B_d(A, \epsilon) = \{x \in X: d(x, A) < \epsilon\} = \cup_{x\in A}B_d(x,\epsilon)$. More generally, for a function $f:(X,d) \to (0,\infty)$, we define the \textit{functional enlargement} (or \textit{pointwise enlargement}) of $A\subseteq X$ by $B_d(A,f) = \cup_{x\in A}B_d(x,f(x))$. It is easy to see that $\overline{A} \subseteq B_d(A, \epsilon)$ for any $\epsilon > 0$, while $\overline{A}$ need not be contained in $B_d(A, f)$ for a general function $f:(X,d) \to (0,\infty)$. Recently in \cite{agarwal2024set}, the authors exhibited that functional enlargements arise naturally in the study of $\tau_{\mathcal{S},d}$-convergence.

	\begin{definition}\label{definition T_S,d}\normalfont (\cite{cao})
		Let $(X,d)$ be a metric space and $\mathcal{S} \subseteq \mathcal{P}_0(X)$. Then a net $(A_\lambda)$ in $\mathcal{P}_0(X)$ is \textit{$\tau_{\mathcal{S},d}^-$-convergent} to a nonempty set $A$ if for any $\epsilon > 0$ and $S \in \mathcal{S}$, eventually, $d(x,A_{\lambda}) - d(x, A) < \epsilon$ for all $x \in S$. The corresponding topology is denoted by  $\tau_{\mathcal{S},d}^-$.
		Similarly, a net $(A_{\lambda})$ in $\mathcal{P}_0(X)$ is \textit{$\tau_{\mathcal{S},d}^+$-convergent} to a nonempty set $A$ if and only if for any $\epsilon > 0$ and $S \in \mathcal{S}$, eventually, $d(x,A) - d(x, A_{\lambda}) < \epsilon$ for all $x \in S$. The corresponding topology is denoted by  $\tau_{\mathcal{S},d}^+$. 
	\end{definition}
	The topology $\tau_{\mathcal{S},d}$ is the join of $\tau_{\mathcal{S},d}^-$ and $\tau_{\mathcal{S},d}^+$. A base for a compatible uniformity for $\tau_{\mathcal{S},d}$ is given by the family $\{U_{\mathcal{S},\epsilon} : S \in \mathcal{S}, \epsilon > 0\}$, where $$ U_{\mathcal{S},\epsilon} = \{(A,B) \in \mathcal{P}_0(X) \times \mathcal{P}_0(X) : |d(x,A) - d(x,B)| < \epsilon\ \ \forall x \in S\}.$$

	
	\begin{definition}\label{S-convergence}\normalfont(\cite{BCL})
		Suppose $\mathcal{S}$ is a nonempty family of subsets of a metric space $(X,d)$. A net $(A_\lambda)$ is said to be \textit{lower bornological convergent} to $A$ in $\mathcal{P}_0(X)$ denoted by $\mathcal{S}^-$-convergence, if for each $\epsilon > 0$ and $S \in \mathcal{S}$, eventually $ A \cap S \subseteq B_d(A_\lambda, \epsilon)$.
		
		We say $(A_\lambda)$ is \textit{upper bornological convergent} to $A$ in $\mathcal{P}_0(X)$ denoted by $\mathcal{S}^+$-convergence, if for each $\epsilon > 0$ and $S \in \mathcal{S}$, eventually $ A_\lambda \cap S \subseteq B_d(A, \epsilon)$.
		
		A net $(A_\lambda)$ $\mathcal{S}$-converges to $A$ provided it is both $\mathcal{S}^-$-convergent and $\mathcal{S}^+$-convergent to $A$. The $\mathcal{S}$-convergence is known as the bornological convergence.
	\end{definition}
	Note that $e_d(A,B) \leq \epsilon \Leftrightarrow A \subseteq B_d(B, \epsilon)$, where $A,B \in \mathcal{P}_0(X)$ and $\epsilon > 0$. So a net $(A_{\lambda})$ is $\mathcal{S}$-convergent to $A$ in $CL(X)$ if and only if 
	for each $S \in \mathcal{S}$, $\lim_{\lambda}e_d(A\cap S, A_{\lambda})  = 0$ and $\lim_{\lambda}e_d(A_{\lambda}\cap S, A)  = 0$. 
	 
	Indeed, for $\mathcal{S} = \mathcal{B}_d(X)$ the $\mathcal{S}$-convergence corresponds to the Attouch-Wets topology $\tau_{AW_d}$. However, in general, the $\mathcal{S}$-convergence may not correspond to a topology. In the last two decades, several researchers have explored structural properties of $\mathcal{S}$ under which $\mathcal{S}$-convergence becomes topological \cite{Bcas,PbciAWc,BCL}.
	
	\begin{definition}\normalfont (\cite{beer2014excess})\label{definition leftexcess}
		Let $(X,d)$ be a metric space and let $\mathcal{S}\subseteq \mathcal{P}_0(X)$. Then the \textit{left excess convergence} on $\mathcal{P}_0(X)$ is the convergence corresponding to the weakest topology on $\mathcal{P}_0(X)$ such that for each $S \in \mathcal{S}$, $$A \mapsto e_d(S, A) \quad(A \in \mathcal{P}_0(X))$$ is a continuous extended real-valued function on $\mathcal{P}_0(X)$. This convergence is denoted by $\mathcal{LE}_{\mathcal{S},d}$. It can be split into two parts:
		\begin{enumerate}[(i)]
			\item \textit{Upper left excess convergence}: It is the convergence corresponding to the weakest topology on $\mathcal{P}_0(X)$ such that each member of the family $\{e_d(S, \cdot): S \in \mathcal{S}\}$ is lower semi-continuous. It is denoted by $\mathcal{LE}_{\mathcal{S},d}^+$. 
			\item \textit{Lower left excess convergence}: It is the convergence corresponding to the weakest topology on $\mathcal{P}_0(X)$ such that each member of the family $\{e_d(S, \cdot): S \in \mathcal{S}\}$ is upper semi-continuous. It is denoted by $\mathcal{LE}_{\mathcal{S},d}^-$. 
		\end{enumerate}
	\end{definition} 

	It is well-known that for a family $\mathcal{S}$ of nonempty subsets of $(X,d)$ containing all singleton subsets of $X$, the $\mathcal{LE}_{\mathcal{S},d}^+$-convergence always coincides with the upper Wijsman convergence. But  $\mathcal{LE}_{\mathcal{S},d}^-$-convergence is sensitive to the choice of family $\mathcal{S}$ (see the end of this section). 
	
	As mentioned in pg $115$ of \cite{ToCCoS}, the full $\mathcal{LE}_{\mathcal{S},d}$-convergence in particular cases such as $\mathcal{S} = \mathcal{B}_d(X)$ and $\mathcal{P}_0(X)$ form natural dual convergences to the bounded proximal and proximal convergences, respectively. Therefore, for $\mathcal{S} = \mathcal{B}_d(X)$ and $\mathcal{P}_0(X)$, $\mathcal{LE}_{\mathcal{S},d}$-convergences are called as \textit{dual bounded proximal convergence} and \textit{dual proximal convergence}, respectively. In fact the join of bounded proximal convergence (proximal convergence) and dual bounded bounded proximal convergence (dual proximal convergence) amounts to the Attouch-Wets convergence (Hausdorff metric convergence) on $CL(X)$. 
	

	In general, we have the following relations among the above defined set convergences (see Proposition $2$ of \cite{gapexcess}). $$ \mathcal{S}^- \leq \mathcal{LE}_{\mathcal{S},d}^- \leq \tau_{\mathcal{S},d}^-; \text{ and }   \mathcal{S}^+ \leq \tau_{\mathcal{S},d}^+.$$ However, $\mathcal{LE}_{\mathcal{S},d}^+$ and $\mathcal{S}^+$ may not be related in general. It is to be noted that many set convergences such as hit-and-miss, proximal hit-and-miss and gap convergences have their lower parts same as well-known lower Vietoris convergence. But this is not true for the convergences under consideration in the present paper. In fact, in the three standard cases when $\mathcal{S} = \mathcal{F}(X), \mathcal{B}_d(X)$, and $\mathcal{P}_0(X)$, we have 
	\begin{itemize}
		\item $\tau_{\mathcal{S},d}^{-} = \mathcal{S}^{-} = \mathcal{LE}_{\mathcal{S},d}^{-} = \tau_{W_d}^{-}$ for $\mathcal{S} = \mathcal{F}(X)$;
		\item $\tau_{\mathcal{S},d}^{-} = \mathcal{S}^{-} = \mathcal{LE}_{\mathcal{S},d}^{-} = \tau_{{AW}_d}^{-}$ for $\mathcal{S} = \mathcal{B}_d(X)$;
		\item $\tau_{\mathcal{S},d}^{-} = \mathcal{S}^{-} = \mathcal{LE}_{\mathcal{S},d}^{-} = \tau_{H_d}^{-}$ for $\mathcal{S} = \mathcal{P}_0(X)$;
	\end{itemize}

	\section{$\tau_{\mathcal{S},d}$-convergence vis-\'a-vis $\mathcal{S}$-convergence}
	The  main aim of this section is to examine the equivalence of $\tau_{\mathcal{S},d}$-convergence and $\mathcal{S}$-convergence. It is shown that $\tau_{\mathcal{S},d} = \mathcal{S}$ if and only if $\tau_{\mathcal{S},d}^+ = \mathcal{S}^+$ and $\tau_{\mathcal{S},d}^- = \mathcal{S}^-$.  In \cite{agarwal2024set}, the authors present a necessary and sufficient condition on the structure of the bornology $\mathcal{S}$ for $\tau_{\mathcal{S},d}^+ = \mathcal{S}^+$ on $CL(X)$. In this section, we provide a characterization for the coincidence  $\tau_{\mathcal{S},d}^- = \mathcal{S}^-$. We also give some stability  conditions on the underlying bornology that enforce the equivalence of the above mentioned convergences. 
	
	We start with the following perspective of $\mathcal{S}^-$-convergence on $CL(X)$.
	\begin{theorem}[Hit-type characterization for $\mathcal{S}^-$-convergence]\label{Hit-type characterization for S-convergence}
		Let $(X,d)$ be a metric space and  let $\mathcal{S}$ be a bornology on $X$. Suppose $(A_{\lambda})$ is a  net in $CL(X)$ and $A \in CL(X)$. Then the following statements are equivalent:
		\begin{enumerate}[(i)]
			\item $(A_{\lambda})$ is $\mathcal{S}^-$-convergent to $A$;
			\item for $S \in \mathcal{S}$ and $\epsilon > 0$, whenever $A \cap S \neq \emptyset$ then eventually $A_{\lambda}\cap B_d(x, \epsilon) \neq \emptyset$ $\forall x \in A\cap S$. 
		\end{enumerate}
	\end{theorem}
	\begin{proof}
		$(i)\Rightarrow(ii).$ Let $S \in \mathcal{S}$ and  $\epsilon > 0$. If $A \cap S = \emptyset$, then we are done.  Otherwise by $(i)$, there is an $\lambda_0$ such that $A \cap S \subseteq B_d(A_{\lambda}, \epsilon)$ $\forall \lambda \geq \lambda_0$. Consequently, for $x \in A\cap S$ and $\lambda \geq \lambda_0$, $d(x, A_{\lambda}) < \epsilon$. Thus $(ii)$ follows. 
		
		$(ii)\Rightarrow (i).$ Let $S \in \mathcal{S}$ and $\epsilon > 0$. If $A \cap S =\emptyset$, then we are done. Otherwise choose $\lambda_0$ such that $\lambda \geq \lambda_0$, we have $A_{\lambda} \cap B_d(x, \epsilon) \neq \emptyset$ $\forall x \in A \cap S$. Pick $x \in A \cap S$. Then for each $\lambda \geq \lambda_0$, there is an $a_{\lambda} \in A_{\lambda}$ such that $d(x, a_{\lambda}) < \epsilon$. Consequently, $A \cap S \subseteq B_d(A_{\lambda}, \epsilon)$ $\forall \lambda \geq \lambda_0$.  
	\end{proof}
	It is evident that Theorem \ref{Hit-type characterization for S-convergence} holds for an arbitrary family $\mathcal{S}\subseteq \mathcal{P}_0(X)$. Moreover, it follows that for any $A \in CL(X)$ the net $(A \cap S)_{A\cap S \neq \emptyset}$ is $\mathcal{S}^-$-convergent to $A$, where $\mathcal{S}$ is directed by set inclusion.
	
	Now, we have the following hit-and-miss type characterization for $\mathcal{S}$-convergence which follows from Theorem $6.1$ of \cite{agarwal2024set} and Theorem \ref{Hit-type characterization for S-convergence}.  
	
	\begin{theorem}[Hit-and-Miss type characterization of $\mathcal{S}$-convergence]
		Let $(X,d)$ be a metric space and let $\mathcal{S}$ be a bornology on $X$. Suppose $(A_{\lambda})$ is a net in $CL(X)$ and $A \in CL(X)$. Then the following statements are equivalent:
		\begin{enumerate}[(i)]
			\item $(A_{\lambda})$ $\mathcal{S}$-converges to $A$;
			\item for $S_1, S_2 \in \mathcal{S}$ and $\alpha$, $\epsilon > 0$, whenever $A \cap B_d(S_1, \alpha) = \emptyset$ and $A \cap S_2 \neq \emptyset$, then eventually $A_{\lambda} \cap S_1 = \emptyset$ and $A_{\lambda} \cap B_d(x, \epsilon) \neq \emptyset$ $\forall x \in A\cap S_2$. 
		\end{enumerate}
	\end{theorem}
	
	
 Recall that for a bornology $\mathcal{S}$ on $(X,d)$, $\mathcal{S}^- \leqslant \mathcal{LE}_{\mathcal{S},d}^-\leqslant \tau_{\mathcal{S},d}^-$ on $CL(X)$.	For the sake of readers' interest, we state some results from the literature regarding the equivalence of  $\mathcal{S}^-$-convergence with $\mathcal{LE}_{\mathcal{S},d}^-$ and  $\tau_{\mathcal{S},d}^-$.   A strong sufficient condition for the coincidence of $\tau_{\mathcal{S},d}^-$ and $\mathcal{S}^-$ is that $\mathcal{S} \subseteq \mathcal{TB}_d(X)$. As $\mathcal{S}^- = \tau_{W_d}^- = \tau_{\mathcal{S},d}^-$ provided $\mathcal{S} \subseteq \mathcal{TB}_d(X)$ (see, Corollary 27 of \cite{Idealtopo}). 
	 
	 The necessary and sufficient conditions for the coincidence of $\mathcal{LE}_{\mathcal{S},d}^-$ and $\tau_{\mathcal{S},d}^-$ with $\mathcal{S}^-$-convergence were explored by Beer et al. in \cite{gapexcess, Idealtopo}. More precisely, they characterized these coincidences in terms of convergence of certain nets as given in the following result.
	  \begin{proposition}\label{beercoincidenceresult} $($\cite{gapexcess, Idealtopo}$)$
	  	Let $(X,d)$ be a metric space and let $\mathcal{S}$ be a bornology on $X$.
	  	 Then the following statements are equivalent:
	  	 \begin{enumerate}[(i)]
	  	 	\item for each $A \in CL(X)$, the net $(A\cap S)_{A\cap S \neq \emptyset}$ is $\tau^-$-convergent to $A$;
	  	 	\item $\tau^- = \mathcal{S}^-$,
	  	 	
	  	 \end{enumerate}
	  	  where $\tau^- = \tau_{\mathcal{S},d}^-$ and $\mathcal{LE}_{\mathcal{S},d}^-$, respectively. 
	  \end{proposition}
	   
	   The condition $(i)$ of Proposition \ref{beercoincidenceresult} does not picture out any structural property of bornology $\mathcal{S}$ for the coincidence of these convergences. However, for a bornology $\mathcal{S}$ on $(X,d)$, the coincidence of $\mathcal{S}^-$ and $\mathcal{LE}_{\mathcal{S},d}^-$ have the following nice characterization.
	   
	   \begin{theorem}$($Theorem 4, \cite{gapexcess}$)$
	   	Let $(X,d)$ be a metric space and let $\mathcal{S}$ be a bornology on $X$.
	   	Then the following statements are equivalent:
	   	\begin{enumerate}[(i)]
	   		\item for each $A \in CL(X)$, $S \in \mathcal{S}$, whenever $0< \alpha< \epsilon$ such that $S \subseteq B_d(A,\alpha)$, then $\exists$ $S' \in \mathcal{S}$ such that $S \subseteq B_d(A\cap S', \epsilon)$;   
	   		\item $\mathcal{LE}_{\mathcal{S},d}^- = \mathcal{S}^-$.
	   		\end{enumerate}
	   \end{theorem}  
	    A comparatively less technical sufficient condition for the coincidence of $\mathcal{LE}_{\mathcal{S},d}^-$ and $\mathcal{S}^-$ is stated in the next result.
	    \begin{proposition}\label{stable under enlargement implies bornological and excess convergenec are same}$($Corollary $3$ of \cite{gapexcess}$)$
	    	Let $(X,d)$ be a metric space and let $\mathcal{S}$ be a bornology on $X$. If $\mathcal{S}$ is stable under enlargements, that is, if for each $S \in \mathcal{S}$ and $\epsilon > 0$ $B_d(S, \epsilon) \in \mathcal{S}$, then $\mathcal{LE}_{\mathcal{S},d}^- = \mathcal{S}^-$.
	    	\end{proposition}  
    	In Theorem \ref{lowerexcessandSselections}, we give a weaker condition than the stability of $\mathcal{S}$ under enlargements for the coincidence $\mathcal{LE}_{\mathcal{S},d}^- = \mathcal{S}^-$.
    	
    	For the coincidence of $\mathcal{S}^-$ and $\tau_{\mathcal{S},d}^-$-convergences, the stability of $\mathcal{S}$ under constant enlargements is not enough. This can be seen through the following example. 
	    
	    \begin{example}\label{counterexampleTsdandSlowerexcess}
	    	Let $(X,d) = (\mathbb{R}^2, d_e)$, where $d_e$ denotes the Euclidean metric. Consider the bornology $\mathcal{B}(\mathbb{R}^2)$ having base $\mathcal{B}_0 = \{S \subseteq \mathbb{R}^2:  \exists m \in \mathbb{N} \text{ such that }S \subseteq [-m,m] \times \mathbb{R}\}$. One can easily verify that $\mathcal{B}(\mathbb{R}^2)$ is stable under all constant enlargements. However, $\tau_{\mathcal{B}(\mathbb{R}^2),d}^- \neq \mathcal{B}(\mathbb{R}^2)^-$ on $CL(X)$. To see this, let $A$ be the set of all points on the line $y = x$ in the first quadrant and for each $n \in \mathbb{N}$, let $A_n = \{(x,x+ \frac{x}{n}): x \in [0,n]\}$. Take $S \in \mathcal{B}_0$ and $\epsilon > 0$. Choose $m, n_0\in \mathbb{N}$ such that $S \subseteq [-m,m]\times \mathbb{R}$ and $n_0 > \frac{m}{\epsilon}$. Pick $(x,x) \in A \cap S$. Then $d_e((x,x), (x,x+\frac{x}{n})) < \epsilon$ $\forall$ $n \geq n_0$. So $A_n \cap B_d((x,x), \epsilon) \neq \emptyset$ $\forall$ $n \geq n_0$. Thus, $(A_n)$ is $\mathcal{B}(\mathbb{R}^2)^-$-convergent to $A$ in $CL(X)$. On the other hand, for $S=\{(0,n): n \in \mathbb{N}\}$, $d_e((0,2n), A_n)- d_e((0,2n), A) > \frac{1}{2}$. Therefore, $(A_n)$ is not $\tau_{\mathcal{B}(\mathbb{R}^2),d}^-$-convergent to $A$. 
	    	\qed 
	    \end{example}
 	  
	In order for further investigations, we first introduce the following terminology. Let $\mathcal{Z}^+ = \{f:(X,d)\to \mathbb{R}: f(x) > 0 \text{ for all } x\in X\}$ and $L(\mathcal{Z}^+) = \{f\in \mathcal{Z}^+: f \text{ is Lipschitz}\}$. 
		\begin{definition}\normalfont
		Let $(X,d)$ be a metric space. Given nonempty subsets $A$ and $C$ of $X$ and $f \in \mathcal{Z}^+$, we call \textit{$A$ hits $f$-enlargement of $C$ pointwise} if $A \cap B_d(x,f(x)) \neq \emptyset$ $\forall x \in C$.    
	\end{definition}
	Note that if $f \in \mathcal{Z}^+$ is a constant function, then $A$ hits $B_d(C,f)$ pointwise provided $C \subseteq B_d(A,f)$. 
	
		\begin{theorem}\label{coincicdence of lower Tsd and S}
		Let $(X,d)$ be a metric space and let $\mathcal{S}$ be a bornology on $X$. Then the following statements are equivalent:
		\begin{enumerate}[(i)]
			\item $\tau_{\mathcal{S},d}^- = \mathcal{S}^-$;
			\item for $S \in \mathcal{S}$, $A \in CL(X)$, and $f,g \in \mathcal{Z}^+$ with $\displaystyle{\inf_{x \in S}(g(x)- f(x))>0}$, whenever $A$ hits $B_d(S,f)$ pointwise, then there exists $S' \in \mathcal{S}$ such that $A \cap S'$ hits $B_d(S,g)$ pointwise;  
			\item for $S \in \mathcal{S}$, $A \in CL(X)$, and $f,g \in L(\mathcal{Z}^+)$ with $\displaystyle{\inf_{x \in S}(g(x)- f(x))>0}$, whenever $A$ hits $B_d(S,f)$ pointwise, then there exists $S' \in \mathcal{S}$ such that $A \cap S'$ hits $B_d(S,g)$ pointwise.   
		\end{enumerate}
	\end{theorem}	 
	
	\begin{proof}
		$(i) \Rightarrow (ii)$. Suppose $(ii)$ fails. Then there exist $S_0 \in \mathcal{S}$, $A \in CL(X)$ and $f,g \in \mathcal{Z}^+$ such that $\inf_{x \in S_0}(g(x)-f(x)) = \epsilon>0$ and $A$ hits $B_d(S_0,f)$ pointwise but for any $S \in \mathcal{S}$, $A \cap S$ does not hit $B_d(S_0,g)$ pointwise. So for any $S \in \mathcal{S}$, we can find an $x_S \in S_0$ satisfying $A \cap S \cap B_d(x_S, g(x_S)) = \emptyset$. Direct $\mathcal{S}$ by set inclusion. For each $S \in \mathcal{S}$, define $A_S = \overline{A \cap S}$. Then $(A_S)_{A\cap S \neq \emptyset}$ is a net in $CL(X)$. We know that the net $(A_S)_{A\cap S \neq \emptyset}$ is $\mathcal{S}^-$-convergent to $A$. We show that it fails $\tau_{\mathcal{S},d}^-$-convergence. To see this, observe that $g(x_S) \leq d(x_S,A_S)$ for any $S \in \mathcal{S}$ and $f(x_S) > d(x_S,A)$. So $d(x_S,A_S) - d(x_S,A) > \epsilon$ for any $S \in \mathcal{S}$. Thus $(A_S)$ is not $\tau_{\mathcal{S},d}^-$-convergent to $A$. 
		
		$(ii)\Rightarrow (iii)$. It is immediate.
		
		$(iii) \Rightarrow (i)$. Let $(A_\lambda)$ be a net in $CL(X)$ which is $\mathcal{S}^-$-convergent to $A \in CL(X)$. Take $S \in \mathcal{S}$ and $\epsilon > 0$. Define $f = d(\cdot,A) + \epsilon$ and $g = d(\cdot, A) + 2 \epsilon$, so $f,g \in L(\mathcal{Z}^+)$ with $\inf_{x \in S}(g(x)-f(x)) = \epsilon$. Clearly $A$ hits $B_d(S,f)$ pointwise. By $(iii)$, there exists $S' \in \mathcal{S}$ such that $A \cap S'$ hits $B_d(S,g)$ pointwise. Also $(A_{\lambda})$ is $\mathcal{S}^-$-convergent to $A$, $\exists$ $\lambda_0$ such that $A\cap S' \subseteq B_d(A_{\lambda}, \epsilon)$ $\forall$ $\lambda \geq \lambda_0$. Take $x \in S$ and $\lambda \geq \lambda_0$. Then we can find a $y_x \in A \cap S'$ such that $d(x,y_x) < g(x)$. Also $d(y_x, A_{\lambda}) < \epsilon$. Then $d(x,A_{\lambda}) \leq d(x,y_x)+d(y_x, A_{\lambda}) < g(x)+ \epsilon$. Consequently, $d(x,A_{\lambda}) < d(x, A) + 3 \epsilon$ $\forall$ $\lambda \geq \lambda_0$ and $x \in S$. Therefore, $(A_{\lambda})$ is $\tau_{\mathcal{S},d}^-$-convergent to $A$.     
	\end{proof} 
	
	\begin{corollary}
		Let $(X,d)$ be a metric space. Then $\tau_{H_d}^- = \mathcal{P}_0(X)^-$. 
	\end{corollary}
	The proof of the following corollary highlights a reason for considering Lipschitz functions ($L(\mathcal{Z}^+)$) in Theorem \ref{coincicdence of lower Tsd and S}.
	
\begin{corollary}\label{coincidence lower BdX Tsd and S} $($see Proposition 3.1.6 of \cite{ToCCoS}$)$
Let $(X,d)$ be a metric space. Then $\tau_{AW_{d}}^- = \mathcal{B}_d(X)^-$ on $CL(X)$. 
\end{corollary}
	
	\begin{proof}
		Let $S \in \mathcal{B}_d(X)$, $A\in CL(X)$ and $f,g \in L(\mathcal{Z}^+)$ such that $\displaystyle{\inf_{x \in S}}\{g(x)-f(x)\}> 0$
		and $A$ hits $B_d(S,f)$ pointwise. For each $x \in S$, choose $a_x \in A\cap B_d(x,f(x))$. Set $S' = \{a_x: x \in S\}$. Note that $A\cap S'$ hits $B_d(S,g)$ pointwise. To finish the proof, we show that $S' \in \mathcal{B}_d(X)$. Take $a_x,a_{x'} \in S'$. Then  $d(a_x, a_{x'}) \leq d(a_x,x)+d(x,x')+d(x',a_{x'}) < f(x)+f(x')+$ diam $(S)$, where diam $(S)$ is diameter of $S$. Since $f \in L(\mathcal{Z}^+)$, by Theorem $9.1$ of \cite{beer2023bornologies}, $f(S)$ is bounded in $\mathbb{R}$. Also diam $(S) < \infty$. So $d(a_x,a_{x'}) < \infty$. Consequently, $S' \in \mathcal{B}_d(X)$. Therefore, by Theorem \ref{coincicdence of lower Tsd and S}, $\tau_{AW_{d}}^- = \mathcal{B}_d(X)^-$.   
		\end{proof}

\begin{corollary}\label{lowerTsd and S}
Let $(X,d)$ be a metric space and let $\mathcal{S}$ be a bornology on $X$. Then the following statements are equivalent:
\begin{enumerate}[(i)]
\item $\mathcal{S}^- \geqslant \tau_{\mathcal{S},d}^-$;				
\item $\mathcal{S}$-convergence $\geqslant \tau_{\mathcal{S},d}^-$-convergence;
\item for $S \in \mathcal{S}$, $A \in CL(X)$, and $f,g \in \mathcal{Z}^+$ with $\displaystyle{\inf_{x \in S}(g(x)- f(x))>0}$, whenever $A$ hits $B_d(S,f)$ pointwise, then there exists $S' \in \mathcal{S}$ such that $  A \cap S'$ hits $B_d(S,g)$ pointwise.
\end{enumerate}
\end{corollary}
		
\begin{proof}
The implication	$(i) \Rightarrow (ii)$ is trivial, and $(iii) \Rightarrow (i)$ follows from Theorem \ref{coincicdence of lower Tsd and S}.
We only need to show $(ii) \Rightarrow (iii)$. Suppose $(iii)$ fails. Observe that the net $(A_S)$ constructed in Theorem \ref{coincicdence of lower Tsd and S} is also $\mathcal{S}^+$-convergent to $A$. Consequently, $(A_S)$ is $\mathcal{S}$-convergent to $A$ while the $\tau_{\mathcal{S},d}^-$-convergence of $(A_S)$ to $A$ fails. Thus, we arrive at a contradiction.
\end{proof}
From the proof of Corollary \ref{coincidence lower BdX Tsd and S}, a property of $\mathcal{B}_d(X)$ can be extracted. To do this, we first define a selection of a set.
		
Consider a nonempty subset $S$ of a metric space $(X,d)$ and $f \in \mathcal{Z}^+$. For each $x \in S$, choose $y_x \in B_d(x, f(x))$. Then $S' = \{y_x:  x\in S\}$  is called an \textit{$f$-selection of $S$}. Moreover, when $f \in L(\mathcal{Z}^+)$, then an $f$-selection of $S$ is also phrased as \textit{Lip-selection of $S$}. 
\begin{definition}\label{stable under lip-selections}\normalfont
If a family $\mathcal{S}$ contains Lip-selections of all of its members, then we say \textit{$\mathcal{S}$ is stable under Lip-selections}.
\end{definition} Note that for each $f \in \mathcal{Z}^+$ and $S$ is always an $f$-selection of itself. 
		
\begin{proposition}\label{selectionsforBdX}
Let $(X,d)$ be a metric space. Then $\mathcal{B}_d(X)$ is stable under Lip-selections.
\end{proposition}
\begin{proof}
Let $S \in \mathcal{B}_d(X)$ and $f \in L(\mathcal{Z}^+)$. Let $S'$ be an $f$-selection of $S$. Take $y_x, y_{x'} \in S'$. Then $d(y_x,y_{x'}) \leq d(y_x,x)+ d(x,x')+d(x',y_{x'})\leq f(x)+f(x)+$ diam $(S)$. Proceeding now as in Corollary \ref{coincidence lower BdX Tsd and S}, we get $S' \in \mathcal{B}_d(X)$.
\end{proof}
If for any $S\in \mathcal{S}$ and any constant function $f \in \mathcal{Z}^+$ every $f$-selection of $S$ belongs to $\mathcal{S}$, then we say \textit{$\mathcal{S}$ is stable under constant selections}. Note that if $\mathcal{S}$ is stable under enlargements, then $\mathcal{S}$ is stable under constant selections.  However, the converse is not true. For example each selection of every member of the bornology $\mathcal{F}(\mathbb{R})$ is finite.

Our next result gives a sufficient condition for the coincidence $\mathcal{S}^-=\mathcal{LE}_{\mathcal{S},d}^-$  which is weaker than the condition of Proposition \ref{stable under enlargement implies bornological and excess convergenec are same}. 
		 
\begin{theorem}\label{lowerexcessandSselections}
Let $(X,d)$ be a metric space and let $\mathcal{S}$ be a bornology on $X$. If $\mathcal{S}$ is stable under constant selections, then $\mathcal{LE}_{\mathcal{S},d}^- = \mathcal{S}^-$ on $CL(X)$. 
\end{theorem} 
\begin{proof}
Let $(A_{\lambda})$ be $\mathcal{S}^-$-convergent to $A$ in $CL(X)$. Let $S \in \mathcal{S}$ and $\epsilon > 0$ be such that $e_d(S,A) < \epsilon$. Then for each $x \in S$, $\exists$ $a_x \in A$ such that $d(x, a_x) < \frac{\epsilon}{3}$. Set $S' =\{a_x: x \in S\}$, so $S'$ is a constant selection of $S$. By hypothesis, $S' \in \mathcal{S}$. Since $(A_{\lambda})$ is $\mathcal{S}^-$-convergent to $A$, by Theorem \ref{Hit-type characterization for S-convergence} $ \exists$ $\lambda_0$ such that for each $\lambda \geq \lambda_0$, we have $A_{\lambda} \cap B_d(a_x, \frac{\epsilon}{3}) \neq \emptyset$ $\forall$ $x \in S$. Fix $x \in S$ and $\lambda \geq \lambda_0$. Then $d(x, A_{\lambda}) \leq d(x,a_x)+d(a_x, A_{\lambda}) < \frac{2\epsilon}{3}$. Consequently, $e_d(S, A_{\lambda}) < \epsilon$. Therefore, $(A_{\lambda})$ is $\mathcal{LE}_{\mathcal{S},d}^-$-convergent to $A$ in $CL(X)$.    
\end{proof}
		 

Example \ref{counterexampleTsdandSlowerexcess} can be used to infer that stability under constant selections need not imply the coincidence $\tau_{\mathcal{S},d}^- = \mathcal{S}^-$ on $CL(X)$.  However, it is true under a stronger stability condition. 
\begin{theorem}\label{lowerTsdandSselection}
Let $(X,d)$ be a metric space and let $\mathcal{S}$ be a bornology on $X$. If $\mathcal{S}$ is stable under Lip-selections, then $\tau_{\mathcal{S},d}^- = \mathcal{S}^-$ on $CL(X)$.
\end{theorem}
\begin{proof}
Let $S \in \mathcal{S}$, $A \in CL(X)$, $f,g \in L(\mathcal{Z}^+)$ such that $\displaystyle{\inf_{x \in S}}(g(x) - f(x)) > 0$ and $A$ hits $B_d(S,f)$ pointwise. For each $x \in S$, choose $a_x \in A\cap B_d(x,f(x))$ and set $S_A = \{a_x: x \in S\}$. Then $S_A$ is a Lip-selection of $S$ and by the hypothesis, $S_A \in \mathcal{S}$. So by Theorem \ref{coincicdence of lower Tsd and S}, $\tau_{\mathcal{S},d}^- = \mathcal{S}^-$ on $CL(X)$.   
\end{proof}
\begin{remark}
Note that $\mathcal{F}(X), \mathcal{B}_d(X)$, and $\mathcal{P}_0(X)$ are all stable under Lip-selections. Thus, the coincidence $\tau_{\mathcal{S},d}^- =\mathcal{S}^-$ in these cases follows from Theorem \ref{lowerTsdandSselection}.  
\end{remark} 
The converse of Theorem \ref{lowerTsdandSselection} may not hold in general.
\begin{example}\label{counterexampleselections}
Let $(X,d)$ be the sequence space $(\ell_2,\parallel \cdot \parallel_2)$. By Corollary 27 of \cite{Idealtopo}, we have $\tau_{\mathcal{K}(X),d}^- = \mathcal{K}(X)^- = \tau_{W_d}^-$ on $CL(X)$. We show that there exists $K\mathcal{K}(X)$ such that some Lip-selection of $K$ does not belong to $\mathcal{K}(X)$.  For each $n \in \mathbb{N}$, take $a_n = (1,\frac{1}{2},\ldots, \frac{1}{n}, 0, \ldots) \in \ell_2$ and  $a = (\frac{1}{n}) \in \ell_2$.  It is easy to verify that $(a_n)$ converges to $a$ in $(\ell_2, \parallel \cdot \parallel_2)$. Then $K = \{a_n: n \in \mathbb{N}\} \cup \{a\}$ is compact in $(\ell_2, \parallel \cdot \parallel_2)$. Consider $f\in L(\mathcal{Z}^+)$ given by $f(x) = \parallel x \parallel_2 + 2$ $\forall x \in X$ and the set $\{e_n: n \in \mathbb{N}\}$ in $\ell_2$, where $e_n$ is the sequence having $n^{th}$-term `$1$' and remaining terms `$0$'. Since $\{e_n: n \in \mathbb{N}\}$ is uniformly discrete, it is not compact in $\ell_2$. To complete the argument, we show that $\{e_n: n \in \mathbb{N}\}$ is an $f$-selection of $K$ in $(\ell_2, \parallel \cdot \parallel_2)$. To see this, note that for each $n \in \mathbb{N}$, we have $$\parallel a_n - e_n \parallel_2 = \parallel (1, 1/2,\ldots, 1-1/n,0,0,\ldots)\parallel_2 \leq \parallel a_n \parallel_2 + 1 < f(a_n).$$  \qed
\end{example}	
\begin{remark}
It is worth mentioning that the sufficient condition in Theorem \ref{lowerTsdandSselection} for the coincidence $\tau_{\mathcal{S},d}^- =  \mathcal{S}^-$ on $CL(X)$ need not be sufficient for the coincidence of upper parts, that is, $\mathcal{S}$ is stable under Lip-selections not necessarily implies that  $\tau_{\mathcal{S},d}^+ = \mathcal{S}^+$ on $CL(X)$. Example \ref{counterexampleselections} provides a counterexample for this. To illustrate this, note that by Corollary $6.15$ of \cite{agarwal2024set},  we have $\tau_{\mathcal{F}(X),d}^+ = \mathcal{F}(X)^+$ if either each proper closed ball of $(X,d)$ is finite or $(X,d)$ is bounded. However, for $X = \ell_2$ this is not true. 
\end{remark}
		

We now see the coincidence of full $\tau_{\mathcal{S},d}$ and $\mathcal{S}$-convergences.

From the proof of Theorem $4.1$ of \cite{agarwal2024set}, it can be deduced that it holds even if the functions are considered from the smaller class $L(\mathcal{Z}^+)$ than that of $\mathcal{Z}^+$, that is, a net $(A_{\lambda})$ is $\tau_{\mathcal{S},d}^+$-convergent to $A$ is equivalent to the fact that for any $S \in \mathcal{S}$ and $f,g \in L(\mathcal{Z}^+)$ such that $\inf_{x \in S}(g(x) -f(x)) > 0$, whenever $A$ misses $B_d(S,g)$, then eventually $(A_{\lambda})$ misses $B_d(S,f)$. Using this, Theorem $6.4$ of \cite{agarwal2024set} can be restated as follows: the coincidence $\tau_{\mathcal{S},d}^+ = \mathcal{S}^+$ holds on $CL(X)$ provided for $f,g \in L(\mathcal{Z}^+)$ such that $\inf_{x \in S}(g(x)-f(x))> 0$ and $S \in \mathcal{S}$, either $B_d(S,g) = X$ or $B_d(S,f) \in \mathcal{S}$. 

The following theorem gives a sufficient condition for the coincidence of full $\tau_{\mathcal{S},d}$ and $\mathcal{S}$-convergences  
		\begin{theorem}\label{suff Tsd and S}
			Let $(X,d)$ be a metric space and let $\mathcal{S}$ be a bornology on $X$. Suppose $A\in CL(X)$, $(A_{\lambda})$ is a net in $CL(X)$. If for each $S \in \mathcal{S}$ and $f \in L(\mathcal{Z}^+)$ $B_d(S,f) \in \mathcal{S}$, then the following statements are equivalent:
			\begin{enumerate}[(i)]
				\item $(A_{\lambda})$ is $\tau_{\mathcal{S},d}$-convergent to $A$;
				\item $(A_{\lambda})$ is $\mathcal{S}$-convergent to $A$;
				\item for each $S \in \mathcal{S}$, $\lim_{\lambda}e_d(A\cap S, A_{\lambda})  = 0$ and $\lim_{\lambda}e_d(A_{\lambda}\cap S, A)  = 0$;
				\item for each $S$ in some cofinal subset $\mathcal{S}_0$ of $\mathcal{S}$, $\lim_{\lambda}e_d(A\cap S, A_{\lambda})  = 0$ and $\lim_{\lambda}e_d(A_{\lambda}\cap S, A)  = 0$.
			\end{enumerate}  
		\end{theorem}
		\begin{proof} Since $\tau_{\mathcal{S},d} \geqslant \mathcal{S}$, we have $(i) \Rightarrow (ii)$.
			
			$(ii) \Rightarrow (i)$.
		 Let $(A_{\lambda})$ be a net that $\mathcal{S}$-converges to $A$ in $CL(X)$. Then $(A_{\lambda})$ is $\mathcal{S}^+$-convergent to $A$. By Theorem $6.4$ of \cite{agarwal2024set} and the above discussion, $(A_{\lambda})$ is $\tau_{\mathcal{S},d}^+$-convergent to $A$. For lower parts, given that for each $S \in \mathcal{S}$ and $f \in L(\mathcal{Z}^+)$, $B_d(S,f) \in \mathcal{S}$, so all Lip-selections of $S$ lie in $\mathcal{S}$. Consequently, by Theorem \ref{lowerTsdandSselection}, $\tau_{\mathcal{S},d}^- = \mathcal{S}^-$ on $CL(X)$.  
			Hence, the net $(A_{\lambda})$ is $\tau_{\mathcal{S},d}$-convergent to $A$. 
			
			$(ii) \Leftrightarrow (iii)$. It follows from the discussion done after Definition \ref{S-convergence}.
			
			$(iii) \Leftrightarrow (iv)$ is easy to verify.  
		\end{proof}
		
		
		Note that a net $(A_{\lambda})$ is $\tau_{\mathcal{S},d}^+$-convergent to $A$ if and only if $(A_{\lambda}\cup A)$ is $\tau_{\mathcal{S},d}$-convergent to $A$ \cite{Idealtopo}.
			\begin{corollary}
			Let $(X,d)$ be a metric space and let $\mathcal{S}$ be a bornology on $X$. If for any $S \in \mathcal{S}$ no Lipschitz enlargement of S is equal to $X$, then the following statements are equivalent:
			\begin{enumerate}[(i)]
				\item for $S \in \mathcal{S}$ and $f \in L(\mathcal{Z}^+)$, $B_d(S,f) \in \mathcal{S}$;
				\item $\mathcal{S}$-convergence is compatible with $\tau_{\mathcal{S},d}$-convergence;
				\item $\tau_{\mathcal{S},d}^+ = \mathcal{S}^+$. 
			\end{enumerate}
		\end{corollary}
		
		\begin{proof}The implications $(i) \Rightarrow (ii)$ follows from Theorem \ref{suff Tsd and S} and $(ii) \Rightarrow (iii)$ is immediate.
			
			$(iii) \Rightarrow (i)$. Let $S \in \mathcal{S}$ and $f \in L(\mathcal{Z}^+)$. By Theorem 6.4 of \cite{agarwal2024set}, either $B_d(S,f) \in \mathcal{S}$ or $B_d(S,f+\epsilon) = X$ for any $\epsilon > 0$. By the assumption, $B_d(S,f+\epsilon) \neq X$, thus we have $B_d(S,f) \in \mathcal{S}$. 
		\end{proof}
		Observe that for any unbounded metric space $(X,d)$ and $\mathcal{S} \subseteq \mathcal{B}_d(X)$, no Lipschitz enlargement of any $S \in \mathcal{S}$ is equal to $X$. In particular, for a normed linear space $(X,\parallel\cdot\parallel)$, to show $\mathcal{B}_d(X)$-convergence is compatible with $\tau_{AW_d}$-convergence, it is enough to prove that $\tau_{{AW}_d}^+ = \mathcal{B}_d(X)^+$.

		We now deduce an another necessary and sufficient condition on the structure of bornology $\mathcal{S}$ for the coincidence $\tau_{\mathcal{S},d} = \mathcal{S}$. Before that, we first see when $\mathcal{S}$-convergence is stronger than $\tau_{\mathcal{S},d}^+$ on $CL(X)$. 
		\begin{theorem}\label{upperTsd and S}        
			Let $(X,d)$ be a metric space and let $\mathcal{S}$ be a bornology on $X$. Then the following statements are equivalent:
			\begin{enumerate}[(i)]                      
				\item $ \mathcal{S}^+$-convergence $\geqslant \tau_{\mathcal{S},d}^+$-convergence;
				\item $ \mathcal{S}$-convergence $\geqslant \tau_{\mathcal{S},d}^+$-convergence;
				\item for $S \in \mathcal{S}$ and $f,g \in \mathcal{Z}^+$ with $\displaystyle{\inf_{x \in S}(g(x) -  f(x)) > 0}$ whenever $B_d(S,g) \neq X$, then $B_d(S,f) \in \mathcal{S}$. 
			\end{enumerate}
		\end{theorem}
		\begin{proof}
			The implication	$(i) \Rightarrow (ii)$ is immediate, and $(iii) \Rightarrow (i)$ follows from Theorem $6.4$ of \cite{agarwal2024set}.
			
			$(ii) \Rightarrow (iii)$. Suppose $(iii)$ fails. Then by following the proof of Theorem 6.4 of \cite{agarwal2024set}, we can construct a net $(A_\lambda)$ in $CL(X)$ such that $(A_\lambda)$ is $\mathcal{S}^+$-convergent to  some $A \in CL(X)$ and $A\subseteq A_\lambda$ for each $\lambda$ but $(A_\lambda)$ is not $\tau_{\mathcal{S},d}^+$-convergent to $A$.  Since $A\subseteq A_\lambda$ for each $\lambda$, $(A_\lambda)$ is $\mathcal{S}^-$-convergent to $A$. Therefore, $(A_\lambda)$  is $\mathcal{S}$-convergent to $A$. Hence we arrive at a contradiction.  	\end{proof}
			
		We now give a necessary and sufficient condition for the coincidence $\tau_{\mathcal{S},d} = \mathcal{S}$ on $CL(X)$.  We would like to mention that in Theorem $34.3$ of \cite{beer2023bornologies}, another necessary and sufficient condition, in terms of convergence of certain nets, for the coincidence $\tau_{\mathcal{S},d} = \mathcal{S}$ has been given.
		
		\begin{theorem}\label{Tsd and S}
			Let $(X,d)$ be a metric space and let $\mathcal{S}$ be a bornology on $X$. Then the following statements are equivalent:
			\begin{enumerate} [(i)]
				\item $\tau_{\mathcal{S},d}^+ = \mathcal{S}^+$, $\tau_{\mathcal{S},d}^- = \mathcal{S}^-$;
				\item $\mathcal{S}$-convergence is compatible with $\tau_{\mathcal{S},d}$;
				
				\item the bornology $\mathcal{S}$ has following characteristics:
				\begin{enumerate}[(a)]
					\item for $S \in \mathcal{S}$ and $f,g\in \mathcal{Z}^+$ with $\displaystyle{\inf_{x \in S}(g(x) - f(x))> 0}$ whenever $B_d(S,g)$  $\neq X$, then $B_d(S,f) \in \mathcal{S}$, and;
					\item for $S \in \mathcal{S}$ and $p,q \in \mathcal{Z}^+$ with $\displaystyle{\inf_{x \in S}(q(x) - p(x))> 0}$ whenever $A \in CL(X)$ satisfying $ A$ hits $B_d(S,p)$ pointwise, then there exists $S' \in \mathcal{S}$ such that $A \cap S'$ hits $B_d(S,g)$ pointwise.  			\end{enumerate} 		\end{enumerate}
		\end{theorem}
		
		\begin{proof}
			The implication $(i) \Rightarrow (ii)$ is straightforward. 
			
			$(ii) \Rightarrow (iii)$. Since $\mathcal{S}$-convergence $\geqslant \tau_{\mathcal{S},d}^+$, by Theorem \ref{upperTsd and S}, $(iii)(a)$ follows. Similarly, $\mathcal{S}$-convergence $\geqslant \tau_{\mathcal{S},d}^-$, by Corollary \ref{upperTsd and S}, $(iii)(b)$ follows.
			
			$(iii)\Rightarrow (i)$. By Theorem \ref{upperTsd and S}, $(iii)(a)$ implies $\tau_{\mathcal{S},d}^+ = \mathcal{S}^+$ and by Theorem \ref{coincicdence of lower Tsd and S}, $(iii)(b)$ implies $\tau_{\mathcal{S},d}^- = \mathcal{S}^-$. Thus, $(i)$ follows. 
		\end{proof}
		Note that the bornologies $\mathcal{B}_d(X)$ and $\mathcal{P}_0(X)$ satisfy Theorem \ref{Tsd and S}$(iii)$. So we have the following known result as a corollary to Theorem \ref{Tsd and S} (see \cite{ToCCoS}).  
		\begin{corollary}
			Let $(X,d)$ be a metric space. Then the following statements hold:
			\begin{enumerate}[(i)]
				\item $\mathcal{B}_d(X)$-convergence is compatible with $\tau_{AW_{d}}$;
				\item $\mathcal{P}_0(X)$-convergence is compatible with $\tau_{H_d}$.
			\end{enumerate}  
		\end{corollary}
%
		
		\section{$\tau_{\mathcal{S},d}$-convergence vis-{\'a}-vis gap and excess convergence}
		In \cite{gapexcess}, Beer and Levi investigated a set convergence, known as gap and excess convergence, which is induced by a family of gap and excess functionals having fixed left argument ranging over a bornology $\mathcal{S}$. They provide some nice conditions for the equivalence of this convergence with $\mathcal{S}$-convergence.  In this section, we study gap and excess convergence in relation to $\tau_{\mathcal{S},d}$-convergence. We first recall the definition of gap and excess convergence. 
		
		\begin{definition}\normalfont (\cite{gapexcess})
			Let $(X,d)$ be a metric space and let $\mathcal{S}$ be a bornology on $X$. Then the gap and excess convergence is the convergence corresponding to the weakest topology on $CL(X)$ such that for each $S \in \mathcal{S}$, the maps 
			$$A\to D_d(S,A) \text{ and } A\to e_d(S, A) \quad (A \in CL(X))$$ are continuous on $CL(X)$. It is denoted by $\tau_{GE}^{\mathcal{S}}$.
			\end{definition}
		 In particular, the weak convergence on $CL(X)$ generated by the family $\{D_d(S, \cdot): S \in \mathcal{S}\}$ of gap functionals  is known as the \textit{gap convergence}, denoted by $\mathsf{G}_{\mathcal{S},d}$. The gap convergence was formally studied by Beer et al. in \cite{beer2013gap}, where they investigated the equivalence of gap convergences corresponding to two different metrics (as well as bornologies). Clearly, $\tau_{GE}^{\mathcal{S}} = \mathsf{G}_{\mathcal{S},d} \vee \mathcal{LE}_{\mathcal{S},d}$, where $\mathcal{LE}_{\mathcal{S},d}$ is the left excess convergence (see, Definition \ref{definition leftexcess}). More precisely, by Proposition $4$ of \cite{gapexcess}, $\tau_{GE}^{\mathcal{S}}$ is the join of upper gap convergence ($\mathsf{G}_{\mathcal{S},d}^+$) and lower left excess convergence ($\mathcal{LE}_{\mathcal{S},d}^-$), that is, $\tau_{GE}^{\mathcal{S}} = \mathsf{G}_{\mathcal{S},d}^+ \vee \mathcal{LE}_{\mathcal{S},d}^-$.
		  
		 	\begin{proposition}\label{gapexcessTsdbasicresult} $($Proposition 2 of \cite{gapexcess}$)$
		 	Let $(X,d)$ be a metric space and let $\mathcal{S}$ be a bornology on $X$. Then the following statements are true:
		 	\begin{enumerate}[(i)]
		 		\item $\tau_{\mathcal{S},d}\geqslant \tau_{GE}^\mathcal{S}$ on $CL(X)$ $(\mathcal{P}_0(X))$;
		 		\item $\tau_{\mathcal{S},d}^+ \geqslant \mathsf{G}_{\mathcal{S},d}^+$ on $CL(X)$ $(\mathcal{P}_0(X))$;
		 		\item $\tau_{\mathcal{S},d}^- \geqslant \mathcal{LE}_{\mathcal{S},d}^-$ on $CL(X)$ $(\mathcal{P}_0(X))$.
		 	\end{enumerate}  
		 \end{proposition}

	 In the special cases, such as $\mathcal{S} = \mathcal{F}(X)$, $\mathcal{B}_d(X)$, and $\mathcal{P}_0(X)$ both $\tau_{\mathcal{S},d}$ and $\tau_{GE}^\mathcal{S}$ reduce to the classical  Wijsman,  Attouch-Wets, and the Hausdorff metric convergences, respectively. However, the coincidence $\tau_{\mathcal{S},d} = \tau_{GE}^\mathcal{S}$ fails for a general bornology $\mathcal{S}$ on $(X,d)$ as shown by the following counterexample. It seems that such an example is not available in the literature. 
	 
	 	\begin{figure} \includegraphics[height= 6.5cm, width = 12.5 cm]{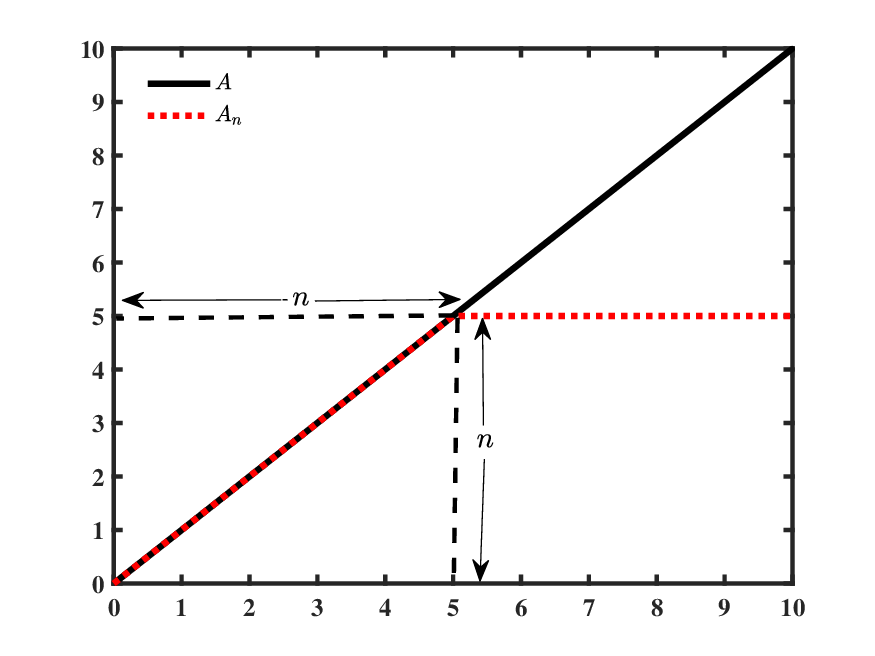} \caption{} \end{figure}
	\begin{example}\label{counterexampleTsdgapexcess}
	Let $X=\{(x,y) \in \mathbb{R}^2: x \geq 0, y \geq 0\}$ equipped with the Euclidean metric. Let $S_0=\{(n,0): n \in \mathbb{N}\}\cup\{(0,n): n \in \mathbb{N}\}$. Consider the bornology $\mathcal{S}$ on $X$ having base as $\{S_0 \cup F: F \in \mathcal{F}(X)\}$. Define $A = \{(x,x) \in \mathbb{R}^2: x \geq 0\}$ and for each $n \in \mathbb{N}$, $A_n = \{(x,x): x \in [0,n]\} \cup \{(x,n): x > n\}$.  Note that by Exercise 2, page 32 in \cite{ToCCoS}, we have $D_d(S_1\cup S_2, C) = \min\{D_d(S_1, C), D_d(S_2, C)\}$ and $e_d(S_1\cup S_2, C) = \max\{e_d(S_1, C), e_d(S_2, C)\}$. So to prove $A_n\xrightarrow{\tau_{GE}^{\mathcal{S}}}A$, it is enough to prove that $$D_d(F,A_n) \to D_d(F,A);  e_d(F,A_n) \to e_d(F,A) ~~  (F \in \mathcal{F}(X)) \text{ and }$$  $$D_d(S,A_n) \to D_d(S,A);  e_d(S,A_n) \to e_d(S,A) ~~ (S\subseteq S_0). \hspace{1.5 cm}$$ 
	
	It is easy to verify that $(A_n)$ is $\tau_{W_d}$-convergent to $A$ in $CL(X)$. So $D_d(F,A) = \lim D_d(F, A_n)$ and $e_d(F,A) = \lim e_d(F, A_n)$ for each $F \in \mathcal{F}(X)$.

	Let $S \subseteq S_0$ and $\alpha >0$. If $D_d(S, A) > \alpha$, then choose $n_0 > 2 \alpha$. So it is clear from Figure $1$ that $D_d(S,A_n) > \alpha$ $\forall$ $n \geq n_0$. Consequently, $D_d(S,A) \leq \liminf D_d(S,A_n)$. Also by Proposition $5$ of \cite{gapexcess} for each $S \in \mathcal{S}$, $D_d(S, \cdot)$ is upper semicontinuous on $(CL(X), \tau_{GE}^{\mathcal{S}})$. So $D_d(S, A)\geq \limsup D_d(S,A_n)$. Thus, $D_d(S, A) = \lim D_d(S, A_n)$. Now if $e_d(S, A) < \alpha$, then $S$ is a bounded subset of $S_0$. Hence $S$ must be finite. So by the above argument, we have $e_d(S, A) = \lim e_d(S, A_n)$. 
	
	We now show that $(A_n)$ is not $\tau_{\mathcal{S},d}$-convergent to $A$ in $CL(X)$. To see this, let $S = \{(n,0): n \in \mathbb{N}\}$ and $f,g \in \mathcal{Z}^+$ be such that $f((x,0)) = \frac{x}{4}$ and $g((x,0)) = \frac{x}{2}$ for $x > 0$. So $A \cap B_d(S,g) = \emptyset$ but for each $n \in \mathbb{N}$, $A_n \cap B_d((25n,0), \frac{25n}{4})) \neq \emptyset$. Therefore, by Theorem $4.1$ of \cite{agarwal2024set},  $\tau_{\mathcal{S},d}^+$-convergence of $(A_n)$ to $A$ fails. \qed  
		\end{example}
	
	 As $\tau_{GE}^\mathcal{S} = \mathsf{G}_{\mathcal{S},d}^+ \vee  \mathcal{LE}_{\mathcal{S},d}^-$, whenever $\tau_{\mathcal{S},d}^+ = \mathsf{G}_{\mathcal{S},d}^+$ and $\tau_{\mathcal{S},d}^- = \mathcal{LE}_{\mathcal{S},d}^-$, then the coincidence $\tau_{\mathcal{S},d}=\tau_{GE}^{\mathcal{S}}$ holds on $CL(X)$. Recently, the authors in \cite{agarwal2024set} furnished necessary and sufficient conditions for the coincidence of $\tau_{\mathcal{S},d}^+$ and $\mathsf{G}_{\mathcal{S},d}^+$. For the other coincidence, we first mention an  important representation for $\mathcal{LE}_{\mathcal{S},d}^-$. It is a consequence of Lemma $3.5$ of \cite{beer2014excess}. 
	
	\begin{theorem}\label{lowerleftexcesschracterization}
		Let $(X,d)$ be a metric space and let $\mathcal{S}$ be a bornology on $X$. Suppose $(A_{\lambda})$ is a net in $CL(X)$ and $A \in CL(X)$. Then the following statements are equivalent:
		\begin{enumerate} [(i)]
			\item $(A_{\lambda})$ is $\mathcal{LE}_{\mathcal{S},d}^-$-convergent to $A$;
			\item whenever $S \in \mathcal{S}$ and $0 <  \alpha < \epsilon$ are such that $A$ hits $B_d(S, \alpha)$ pointwise, then eventually $(A_{\lambda})$ hit $B_d(S, \epsilon)$ pointwise. 
		\end{enumerate}
	\end{theorem} 
%

	Theorem \ref{lowerleftexcesschracterization} gives a hit-type convergence criterion of $\mathcal{LE}_{\mathcal{S},d}^-$. Here we see that the constant enlargements of members of $\mathcal{S}$ play the central role. Surprisingly, such enlargements are not enough for $\tau_{\mathcal{S},d}^-$-convergence. 
	\begin{example}\label{Tsd_excess_counterexample}
		Consider $(X,d)$, $\mathcal{S}$, $A$, and $(A_n)$ as in Example \ref{counterexampleTsdgapexcess}. We show that $(A_n)$ is not $\tau_{\mathcal{S},d}^-$-convergent to $A$. Take $S = \{(0,n): n \in \mathbb{N}\}$. Then for each $n \in \mathbb{N}$, we have $d_e((0,8n), A_n) - d_e((0,8n), A)> 1$. \qed
	\end{example}
	
	If we consider functional enlargements of members of $\mathcal{S}$ in place of their constant enlargements in $(ii)$ of Theorem \ref{lowerleftexcesschracterization} we obtain a new representation for $\tau_{\mathcal{S},d}^-$-convergence. Unlike the usual function space definition of $\tau_{\mathcal{S},d}^-$-convergence, this representation laid emphasis on the hitting behavior of an $\tau_{\mathcal{S},d}^-$-convergent net with functional enlargements of members of $\mathcal{S}$.
		\begin{theorem}[Hit-type characterization for $\tau_{\mathcal{S},d}^-$]\label{Hittypetheoremforlowertsd}
		Let $(X,d)$ be a metric space and let $\mathcal{S}$ be a bornology  on $X$. Suppose $(A_{\lambda})$ is a net in $CL(X)$ and $A \in CL(X)$. Then the following statements  are equivalent:
		\begin{enumerate}[(i)]
			
			\item  $(A_{\lambda})$ $\tau_{\mathcal{S},d}^-$-converges to $A$;
			\item for $S \in \mathcal{S}$ and $f,g \in \mathcal{Z}^+$ with $\displaystyle{\inf_{x \in S}(g(x)-f(x))> 0}$, whenever $A$ hits $B_d(S, f)$ pointwise, then eventually $(A_{\lambda})$ hit $B_d(S, g)$ pointwise; 
			\item for $S \in \mathcal{S}$ and $f,g \in L(\mathcal{Z}^+)$ with $\displaystyle{\inf_{x \in S}(g(x)-f(x))> 0}$, whenever $A$ hits $B_d(S, f)$ pointwise, then eventually $(A_{\lambda})$ hit $B_d(S, g)$ pointwise.
		\end{enumerate}
	\end{theorem}
	\begin{proof}
		$(i)\Rightarrow (ii)$. Let $S \in \mathcal{S}$ and $f,g \in \mathcal{Z}^+$ be such that $\displaystyle{\inf_{x \in S}(g(x) - f(x))}= r > 0$ and $A \cap B_d(x, f(x)) \neq \emptyset$ $\forall x \in S$. By $(i)$, there is a $\lambda_0$ such that $\forall \lambda \geq \lambda_0$, we have $d(x, A_{\lambda}) < d(x, A) + \frac{r}{2}$ $\forall x \in S$. Fix $\lambda \geq \lambda_0$ and $x \in S$. Then $d(x, A_{\lambda}) < d(x, A) + \frac{r}{2} \leq f(x) + \frac{r}{2}$. So $d(x, A_{\lambda}) < g(x)$ $\forall x \in S$. Consequently, for all $\lambda \geq \lambda_0$, we have $A_{\lambda} \cap B_d(x,g(x)) \neq \emptyset$ $ \forall x \in S$. 
		
		$(iii) \Rightarrow (i)$. Let $ S \in \mathcal{S}$ and $ \epsilon > 0$. Define $f(x) = d(x, A) + \frac{\epsilon}{3}$ and $g(x) = d(x, A) + \frac{\epsilon}{2}$ for every $x \in S$. Clearly, $f,g \in L(\mathcal{Z}^+)$ and $A \cap B_d(x, f(x)) \neq \emptyset$ $\forall x \in S$. Then by the hypothesis, $A_{\lambda} \cap B_d(x, g(x)) \neq \emptyset$ $\forall x \in S$ eventually. Consequently, $d(x, A_{\lambda}) \leq g(x) <  d(x, A) + \epsilon$ $\forall x \in S$. Therefore, eventually $\forall x \in S$, we have $d(x, A_{\lambda}) - d(x, A) < \epsilon$. 
	\end{proof}
	\begin{remark}
		Note that the previous theorem holds for an arbitrary family $\mathcal{S} \subseteq \mathcal{P}_0(X)$.
	\end{remark}

	 We now come to one of the main results of our section. Before that we fix some preliminary notations. For a nonempty subset $A$ of $(X,d)$ and $f \in \mathcal{Z}^+$, define $H_d^-(A,f) = \{C \in CL(X): C \text{ hits } B_d(A,f)\text{ pointwise}\}$. When $f$ is a constant function say, $f(x) = \epsilon$ $\forall x$, then $H_d^-(A, \epsilon) = \{C \in CL(X): e_d(A, C) \leq \epsilon\}$. Now Theorem \ref{Hittypetheoremforlowertsd} can be reformulated as: a net $(A_{\lambda})$ is $\tau_{\mathcal{S},d}^-$-convergent to $A$ in $CL(X)$ if and only if whenever $S \in \mathcal{S}$ and $f,g \in \mathcal{Z}^+$ such that $\displaystyle{\inf_{x \in S}}(g(x)-f(x))>0$ and $A \in H_d^-(S,f)$, then $A_{\lambda} \in H_d^-(S,g)$ eventually.  
	 
	 \begin{theorem}\label{lower Tsd and excess}
	 	Let $(X,d)$ be a metric space and let $\mathcal{S}$ be a bornology on $X$. Then the following statements are equivalent. 
	 	\begin{enumerate}[(i)]
	 		\item for  $S \in \mathcal{S}$ and $f,g \in \mathcal{Z}^+$ with $\displaystyle{\inf_{x \in S}(g(x) - f(x))>0}$, whenever $A \in H_d^-(S,f)$ then there exist $0 < \lambda_i < \sigma_i$ and $S_i \in \mathcal{S}$, $i = 1, \ldots, n$ such that $$ A \in \bigcap_{i=1}^{n}H_d^-(S_i, \lambda_i) \subseteq \bigcap_{i=1}^{n}H_d^-(S_i, \sigma_i) \subseteq H_d^-(S,g);$$
	 		
	 		\item for  $S \in \mathcal{S}$ and $f,g \in L(\mathcal{Z}^+)$ with $\displaystyle{\inf_{x \in S}(g(x) - f(x))>0}$, whenever $A \in H_d^-(S,f)$ then there exist $0 < \lambda_i < \sigma_i$ and
	 		$S_i \in \mathcal{S}$, $i = 1, \ldots, n$ such that $$ A \in \bigcap_{i=1}^{n}H_d^-(S_i, \lambda_i) \subseteq \bigcap_{i=1}^{n}H_d^-(S_i, \sigma_i) \subseteq H_d^-(S,g);$$
	 		
	 		\item $\tau_{\mathcal{S},d}^- = \mathcal{LE}_{\mathcal{S},d}^-$ on $CL(X)$.
	 		
	 		\item $\tau_{\mathcal{S},d}^-$ is the weakest topology on $CL(X)$ such that each member of the family of excess functionals $\{e_d(S, \cdot): S \in \mathcal{S}\}$ is upper semicontinuous.  
	 	\end{enumerate}
	 \end{theorem}
	 \begin{proof} It is enough to show $(ii) \Rightarrow (iii)$ and $(iii) \Rightarrow (i)$.

	 	$(ii)\Rightarrow (iii)$. 
	 	Let $(A_{\lambda})$ be a net that $\mathcal{LE}_{\mathcal{S},d}^-$-convergent to $A$ in $CL(X)$. Let $S \in \mathcal{S}$ and $f,g \in L(\mathcal{Z}^+)$ such that $\displaystyle{\inf_{x \in S}(g(x)-f(x) )> 0}$ and $A \in H_d^-(S,f)$. By the hypothesis, there exist $0 < \lambda_i < \sigma_i$ and $S_i \in \mathcal{S}$, $i = 1, \ldots, n$ such that $$ A \in \bigcap_{i=1}^{n}H_d^-(S_i, \lambda_i) \subseteq \bigcap_{i=1}^{n}H_d^-(S_i, \sigma_i) \subseteq H_d^-(S,g).$$ 
	 	By above inclusions, $e_d(S_i, A) < \sigma_i$ $\forall i = 1, \ldots, n$. Since $(A_{\lambda})$ is $\mathcal{LE}_{\mathcal{S},d}^-$-convergent to $A$, $e_d(S_i, A_{\lambda}) < \sigma_i$ eventually $\forall 1\leq i \leq n$. Consequently, $A_{\lambda} \in H_d^-(S_i, \lambda_i)$ eventually for each $i = 1, \ldots, n$. So $A_{\lambda} \in H_d^-(S,g)$ eventually. Hence, $(A_{\lambda})$ is $\tau_{\mathcal{S},d}^-$-convergent to $A$.

	 	$(iii)\Rightarrow (i)$. Suppose $(i)$ fails. Then there exist $S_0 \in \mathcal{S}$, $f,g \in \mathcal{Z}^+$ with $\displaystyle{\inf_{x \in S}(g(x)- f(x))>0}$, and $A \in CL(X)$ such that $A \in H_d^-(S_0,f)$ but for every $0 < \lambda_i < \sigma_i$ and $S_i$, $i = 1, \ldots, n$ satisfying $\displaystyle{A \in  \bigcap_{i=1}^nH_d^-(S_i, \lambda_i)}$ we have $\displaystyle{\bigcap_{i=1}^nH_d^-(S_i, \sigma_i)\nsubseteq H_d^-(S_0,g)}$. 
	 	Set $\mathcal{A} = \{S \in \mathcal{S}: \exists \text{ } \epsilon > 0 \text{ such that } S \subseteq B_d(A, \epsilon)\}$. For $S \in \mathcal{A}$, let $\alpha_S = \min\{\epsilon > 0: S \subseteq B_d(A, \epsilon)\}$. 
	 	Define $$ \Omega = \left\{(F, m) \in \mathcal{F}(\mathcal{A}) \times \mathbb{N}: A \in \bigcap_{S \in F}H_d^-\left(S, \alpha_S + \frac{1}{m}\right) \right\}.$$ Clearly for any $a \in A$ and $m \in \mathbb{N}$, $(\{a\}, m) \in \Omega$, so $\Omega$ is nonempty. Direct $\Omega$ by the relation $\geqslant$: $(F', m') \geqslant (F, m)$ if and only if $F \subseteq F'$ and $m \leq m'$.
	 	By the definition of $\Omega$ and $\alpha_S$, where $S \in \mathcal{A}$, $(\Omega, \geqslant)$ is a directed set. By the assumption, for each $(F, m) \in \Omega$, $\exists$ a $C_{(F, m)} \in \displaystyle{\bigcap_{S \in F}}H_d^-(S, \alpha_S+ 2/m) \setminus H_d^-(S_0,g)$. Consequently, $(C_{(F, m)})_{(F, m) \in \Omega}$ is a net in $CL(X)$. We show that $(C_{(F, m)})$ is $\mathcal{LE}_{\mathcal{S},d}^-$-convergent to $A$ while $\tau_{\mathcal{S},d}^-$-convergence fails. To see $\mathcal{LE}_{\mathcal{S},d}^-$-convergence, let $S_1 \in \mathcal{S}$ and $\epsilon_1 > 0$ be such that $e_d(S_1,A) < \epsilon_1$. 
	 	Choose $m_1 \in \mathbb{N}$ such that $\frac{2}{m_1} < \epsilon_1 - \alpha_{S_1}$. Then for each $(F,m) \in \Omega$ with $(F,m) \geqslant (\{S_1\}, m_1)$, we  have $C_{(F,m)} \in H_d^-\left(S_1, \alpha_{S_1}+\frac{2}{m_1}\right)$. 
	 	So $e_d(S_1, C_{(F,m)}) \leq \alpha_{S_1} + \frac{2}{m_1} < \epsilon_1$ for all $(F,m) \geqslant (\{S_1\}, m)$. 
	 	
	 	For $\tau_{\mathcal{S},d}^-$-convergence, observe that $A \in H_d^-(S_0,f)$  while $C_{(F,m)} \notin H_d^-(S_0,g)$ for any $(F,m) \in \Omega$. Therefore, by Theorem \ref{Hittypetheoremforlowertsd}, the $\tau_{\mathcal{S},d}^-$-convergence of $(C_{(F,m)})$ to $A$ fails. 		
	 \end{proof}
	
		We deduce the following results of Beer and Lucchetti \cite{beer1993weak} as corollaries to Theorem \ref{lower Tsd and excess}.
	 \begin{corollary}\label{lowerattouchexcess}
	 	Let $(X,d)$ be a metric space. Then $\tau_{AW_{d}}^-$ on $CL(X)$ is the weakest topology such that each member of the family of excess functionals $\{e_d(S, \cdot): S \in \mathcal{B}_d(X)\}$ is upper semicontinuous.   
	 \end{corollary}
	 \begin{proof}
	 	It is enough to show that the condition $(ii)$ of Theorem \ref{lower Tsd and excess} holds for $\mathcal{B}_d(X)$.  
	 	Let $S \in \mathcal{B}_d(X)$, $f,g  \in L(\mathcal{Z}^+)$ with $\displaystyle{\inf_{x\in S}}(g(x)-f(x))>0$, and $A \in CL(X)$ such that $A \in H_d^-(S,f)$. Then choose a $f$-selection of $\mathcal{S}$, say $S_A$ such that $S_A \subseteq A$. By Proposition \ref{selectionsforBdX}, $S_A \in \mathcal{B}_d(X)$. Set $2r = \displaystyle{\inf_{x \in S}}(g(x)-f(x))$. So $A \in H_d^-(S_A,r)$. We show that $H_d^-(S_A,2r) \subseteq H_d^-(S,g)$. Let $C \in H_d^-(S_A,2r)$ and $x \in S$. Then $\exists$ $x' \in S_A$ such that $d(x,x')< f(x)$. Consequently, $d(x,C)\leq d(x,x')+ d(x', C)< f(x)+ 2r < g(x)$. So $C \in H_d^-(S,g)$.
	 	\end{proof}

	 		\begin{corollary}
	 		Let $(X,d)$ be a metric space. Then $\tau_{H_d}^-$ on $CL(X)$ is the weakest topology such that each member of the family of excess functionals $\{e_d(S, \cdot): S \in \mathcal{P}_0(X)\}$ is upper semicontinuous. 
	 	\end{corollary}
	 	\begin{proof}
	 		Let $S \in \mathcal{B}_d(X)$, $f,g  \in L(\mathcal{Z}^+)$ with $\displaystyle{\inf_{x\in S}}(g(x)-f(x))>0$, and $A \in CL(X)$ such that $A \in H_d^-(S,f)$. Take $S_A = A \cap B_d(S,f)$. Clearly $S_A \in \mathcal{P}_0(X)$. The rest of the proof is similar to the proof of Corollary \ref{lowerattouchexcess}.
	 	\end{proof}
	 	
	 		Beer et al. in \cite{beer2013gap}, introduced a covering condition on a bornology to study the coincidence  $\mathsf{G}_{\mathcal{S},d}^+ = \mathsf{G}_{\mathcal{S},\rho}^+$ for two metrics $d, \rho$ on $X$, and the coincidence $\mathsf{G}_{\mathcal{S},d}^+ = \mathsf{G}_{\mathcal{T},d}^+$  for two families $\mathcal{S}, \mathcal{T}$ of subsets of $X$. Recently in \cite{agarwal2024set}, the authors used this notion to study the coincidence $\tau_{\mathcal{S},d}^+ = \mathsf{G}_{\mathcal{S},d}^+$.
	 	
	 	
	 	\begin{definition}\normalfont (\cite{beer2013gap})\label{Strictlyincluded}
	 		Let $(X,d)$ be a metric space and let $\mathcal{S}$ be an arbitrary family of nonempty subsets of $X$. We say a subset $A$ of $X$ is \textit{strictly $(\mathcal{S}-d)$ included} in another nonempty subset $C$ of $X$ if there exists a finite subset $\{S_1,\ldots,S_n\}$ of $\mathcal{S}$, and for every $i \in \{1,\ldots,n\}$ there are $\alpha_i, \epsilon_i$ with $0 < \alpha_i < \epsilon_i$, such that $$A \subseteq \cup_{i=1}^nB_d(S_i, \alpha_i)  \subseteq \cup_{i=1}^nB_d(S_i, \epsilon_i) \subseteq C. $$
	 	\end{definition}
	 	
 	From the proof of Theorem $5.5$ of \cite{agarwal2024set} it can be derived that the coincidence $\tau_{\mathcal{S},d}^+ = \mathsf{G}_{\mathcal{S},d}^+$ holds on $CL(X)$ provided for $S\in \mathcal{S}$ and $f,g \in L(\mathcal{Z}^+)$ such that $\inf_{x \in S}(g(x) - f(x))>0$ either $B_d(S,g) = X$ or $B_d(S,f)$ is strictly $(\mathcal{S}-d)$ included in $B_d(S,g)$. 
 	
 	The following theorem provides a two-sided result for the coincidences $\tau_{\mathcal{S},d}^+  = \mathsf{G}_{\mathcal{S},d}^+$ , $\tau_{\mathcal{S},d}^-  = \mathcal{LE}_{\mathcal{S},d}^-$ which follows from Theorem $5.5$ of \cite{agarwal2024set} and Theorem \ref{lower Tsd and excess}.
	 	\begin{theorem}
	 		Let $(X,d)$ be a metric space and let $\mathcal{S}$ be a bornology on $X$. Then the following statements are equivalent:
	 		\begin{enumerate}[(i)]
	 			\item $\tau_{\mathcal{S},d}^+  = \mathsf{G}_{\mathcal{S},d}^+$ , $\tau_{\mathcal{S},d}^-  = \mathcal{LE}_{\mathcal{S},d}^-$; 
	 			\item the bornology $\mathcal{S}$ has the following characteristics:
	 			\begin{enumerate}[(a)]
	 				\item for each $S \in \mathcal{S}$ and $f, g \in L(\mathcal{Z}^+)$ with $\inf\{g(x) - f(x) : x \in S\} > 0$, either $B_d(S , g) = X$ or $B_d(S,f)$ is strictly $(\mathcal{S}-d)$-included in $B_d(S,g)$;
	 				\item for each $S \in \mathcal{S}$ and $f,g \in \mathcal{Z}^+$ with $\displaystyle{\inf_{x \in S}(g(x) - f(x))>0}$, whenever $A \in H_d^-(S,f)$ then there exist $0 < \lambda_i < \sigma_i$ and $S_i \in \mathcal{S}$, $i = 1, \ldots, n$ such that $$ A \in \bigcap_{i=1}^{n}H_d^-(S_i, \lambda_i) \subseteq \bigcap_{i=1}^{n}H_d^-(S_i, \sigma_i) \subseteq H_d^-(S,g).$$
	 			\end{enumerate}
	 			
	 		\end{enumerate}
 		 
	 	\end{theorem}
	 	A less technical sufficient condition for the coincidence $\tau_{\mathcal{S},d}^- = \mathcal{LE}_{\mathcal{S},d}^-$ on $CL(X)$ is given by the following result. 
	 	
	 	\begin{proposition}\label{suffconditionlowerTsdandleftexcess}
	 		Let $(X,d)$ be a metric space and let $\mathcal{S}$ be a bornology on $X$. If $\mathcal{S}$ is stable under Lip-selections, then $\tau_{\mathcal{S},d}^- = \mathcal{LE}_{\mathcal{S},d}^-$ on $CL(X)$. 
	 	\end{proposition}
	 	\begin{proof}
	 		Since $\mathcal{LE}_{\mathcal{S},d}^-$-convergence is stronger than $\mathcal{S}^-$-convergence, Theorem \ref{lowerTsdandSselection} and Proposition \ref{gapexcessTsdbasicresult} together imply that $\tau_{\mathcal{S},d}^- = \mathcal{LE}_{\mathcal{S},d}^- = \mathcal{S}^-$ on $CL(X)$.	 	\end{proof}
	 	
	 	  The stability condition in Proposition \ref{suffconditionlowerTsdandleftexcess} need not be necessary (see, Example \ref{counterexampleselections}).	 
	 	
	 	\begin{remark}\begin{enumerate}[(i)]
	 			\item From Corollary \ref{lowerattouchexcess} and Theorem \ref{lowerleftexcesschracterization}, a net $(A_{\lambda})$ is $\tau_{{AW}_d}^-$-convergent to $A$ provided for each $S \in \mathcal{B}_d(X)$ and $0< \alpha< \epsilon $ whenever $A$ hits $B_d(S, \alpha)$ pointwise, then eventually $(A_{\lambda})$ hits $B_d(S, \epsilon)$ pointwise.
	 			
	 			\item More generally, for a bornology $\mathcal{S}$ on $(X,d)$ which is stable under Lip-selections, it is enough to consider  constant functions in $\mathcal{Z}^+$ to imply $\tau_{\mathcal{S},d}^-$-convergence in Theorem \ref{Hittypetheoremforlowertsd}.  
	 		\end{enumerate}	
	 	\end{remark}
 	
	 	To see the coincidence of $\tau_{\mathcal{S},d}$ and $\tau_{GE}^{\mathcal{S}}$, we need a stronger stability condition on $\mathcal{S}$ than the one given in Proposition \ref{suffconditionlowerTsdandleftexcess}. 
	 	\begin{theorem}\label{Tsdandgapeaxcess}
	 		Let $(X,d)$ be a metric space and let $\mathcal{S}$ be a bornology on $X$. If $\mathcal{S}$ is stable under Lipschitz enlargements, then $\tau_{\mathcal{S},d}^+ = \mathsf{G}_{\mathcal{S},d}^+$, $\tau_{\mathcal{S},d}^- = \mathcal{LE}_{\mathcal{S},d}^-$.
	 	\end{theorem}
	 		\begin{proof}
	 		Since 	
	 		 Since  $\mathcal{S}$ is stable under Lipschitz enlargements, it is stable under enlargements. So by Theorem 5 and Corollary 3 of \cite{gapexcess}, we have $\mathcal{S}^+ = \mathsf{G}_{\mathcal{S},d}^+$, $\mathcal{S}^- = \mathcal{LE}_{\mathcal{S},d}^-$. But by Theorems \ref{suff Tsd and S} and \ref{Tsd and S}, we have  $\tau_{\mathcal{S},d}^+ = \mathcal{S}^+$, and $\tau_{\mathcal{S},d}^- = \mathcal{S}^-$.
	 		\end{proof} 
	 		 \begin{corollary} $($\cite{beer1993weak}$)$
	 			Let $(X,d)$ be a metric space.  Then $\tau_{AW_d}$ is the weak topology generated by the functionals $\{D_d(S, \cdot) : S \in \mathcal{B}_d(X)\} \cup \{e_d(S, \cdot): S \in \mathcal{B}_d(X)\}$.
	 		\end{corollary}
	 		\begin{proof}
	 			By Theorem $9.1$ of \cite{beer2023bornologies}, any Lipschitz function on $S \in \mathcal{B}_d(X)$ is bounded. Then for any $S \in \mathcal{B}_d(X)$ and $f \in L(\mathcal{Z}^+)$, $B_d(S,f) \in \mathcal{B}_d(X)$. So by Theorem \ref{Tsdandgapeaxcess}, we have $\tau_{AW_{d}} = \tau_{GE}^{\mathcal{B}_d(X)}$.
	 		\end{proof}
	 		\begin{corollary} $($\cite{beer1993weak}$)$
	 			Let $(X,d)$ be a metric space.  Then $\tau_{H_d}$ is the weak topology generated by the functionals $\{D_d(S, \cdot) : S \in \mathcal{P}_0(X)\} \cup \{e_d(S, \cdot): S \in \mathcal{P}_0(X)\}$.
	 		\end{corollary}
	 	The next result shows that the covering condition on $\mathcal{S}$ in Definition \ref{Strictlyincluded} is necessary for the coincidence $\tau_{\mathcal{S},d} = \tau_{GE}^{\mathcal{S}}$.

	 	\begin{proposition}\label{TsdGapexcessnecessary}
	 		Let $(X,d)$ be a metric space and let $\mathcal{S}$ be a bornology on $X$, consider the following statements: 
	 		\begin{enumerate}[(i)]
	 			 \item $\tau_{\mathcal{S},d} = \tau_{GE}^{\mathcal{S}}$ on $CL(X)$;
	 			\item for $S \in \mathcal{S}$ and $f,g \in L(\mathcal{Z}^+)$ such that $\inf_{x \in S}(g(x)-f(x))> 0$, either $B_d(S,f)$ is strictly $(\mathcal{S}-d)$ included in $B_d(S,g)$ or $B_d(S,g) = X$.
	 		\end{enumerate}
 		Then $(i) \Rightarrow (ii)$.
	 	\end{proposition}
	 	\begin{proof} 
	 		
	 	 By Theorem 5.5 of \cite{agarwal2024set}, it is enough to show that whenever $\tau_{\mathcal{S},d} = \tau_{GE}^{\mathcal{S}}$ on $CL(X)$, then $\tau_{\mathcal{S},d}^+ = \mathsf{G}_{\mathcal{S},d}^+$ on $CL(X)$. Let $(A_{\lambda})$ be a net in $CL(X)$ which is $\mathsf{G}_{\mathcal{S},d}^+$-convergent to $A \in CL(X)$. Then $(A_{\lambda}\cup A)$ is  $\tau_{GE}^{\mathcal{S}}$-convergent to $A$. Then by $(ii)$, $(A_{\lambda} \cup A)$ is $\tau_{\mathcal{S},d}$-convergent to $A$. Thus, by definition of $\tau_{\mathcal{S},d}^+$-convergence, $(A_{\lambda})$ is $\tau_{\mathcal{S},d}^+$-convergent to $A$.  
	 	\end{proof}
 	
 	\begin{remark}
 		The covering condition on a bornology $\mathcal{S}$ given in Definition \ref{Strictlyincluded} is generally weaker than $\mathcal{S}$ being stable under Lipschitz enlargements. However, if $\mathcal{S}$ is stable under enlargements, then these conditions are easily seen to be equivalent. 	 	\end{remark}
 	
 	We end this paper by giving a covering criterion that provides a sufficient condition for the coincidence  $\tau_{\mathcal{S},d} = \tau_{GE}^{\mathcal{S}}$ on $CL(X)$. This criterion is a stronger form of the covering condition of Definition \ref{Strictlyincluded}.
	 \begin{theorem}
	 Let $(X,d)$ be a metric space and let $\mathcal{S}$ be a bornology on $X$. Suppose for any $S \in \mathcal{S}$, no Lipschitz enlargement of $S$ is equal to $X$. Consider the following statements:
	 	\begin{enumerate}[(i)]
	 		\item for $S \in \mathcal{S}$ and $f,g \in L(\mathcal{Z}^+)$ with $\inf_{x \in S}(g(x)-f(x))> 0$, $\exists$ $S_1, \ldots, S_n \in \mathcal{S}$ and $0 < \alpha_i < \epsilon_i$, $i = 1, \ldots, n$ such that $\max_{1 \leq i \leq n}2 \epsilon_i< \inf_{x \in S}(g(x)-f(x))> 0$ and $$B_d(S,f) \subseteq \cup_{i=1}^nB_d(S_i, \alpha_i) \subseteq \cup_{i=1}^nB_d(S_i, \epsilon_i) \subseteq B_d(S,g);$$
	 		\item $\tau_{\mathcal{S},d}^+ = \mathsf{G}_{\mathcal{S},d}^+$, $\tau_{\mathcal{S},d}^- = \mathcal{LE}_{\mathcal{S},d}^-$ on $CL(X)$;
	 		\item $\tau_{\mathcal{S},d} = \tau_{GE}^{\mathcal{S}}$ on $CL(X)$.
	 	\end{enumerate}
	 	Then $(i) \Rightarrow (ii) \Rightarrow (iii)$ holds.
	 \end{theorem}	
	 \begin{proof}
	 By Proposition \ref{TsdGapexcessnecessary} and the assumption, we have $\tau_{\mathcal{S},d}^+ = \mathsf{G}_{\mathcal{S},d}^+$. It remains to show that $\tau_{\mathcal{S},d}^- = \mathcal{LE}_{\mathcal{S},d}^-$. Let $(A_{\lambda})$ be a net in $CL(X)$ such that $A_{\lambda} \xrightarrow{\mathcal{LE}_{\mathcal{S},d}^-} A \in CL(X)$. Let $S \in \mathcal{S}$ and $f,g \in L(\mathcal{Z}^+)$ such that $\inf_{x \in S}(g(x)-f(x)) = 2 \gamma> 0$ and $A$ hits $B_d(S,f)$ pointwise.  By $(i)$, choose $S_1, \ldots, S_n \in \mathcal{S}$ and $0 < \alpha_i < \epsilon_i$, where $1 \leq i \leq n$ such that $\max_{1\leq i \leq n}\epsilon_i < \gamma$ and $$B_d(S,f) \subseteq \cup_{i = 1}^nB_d(S_i, \alpha_i) \subseteq \cup_{i = 1}^nB_d(S_i, \epsilon_i) \subseteq B_d(S,g).$$  For each $x \in S$, choose $a_x \in A \cap B_d(x, f(x))$.  Take $S' = \cup_{i=1}^nS_i$, $\alpha = \max_{1 \leq i \leq n}\alpha_i$ and $\epsilon = \max_{1 \leq i \leq n}\epsilon_i$. Then for each $x \in S$, there exists $z_x \in S'$ such that $d(a_x, z_x)< \alpha$. Set $S'' = \{z_x: x \in S\}$. So $A$ hits $B_d(S'',\alpha)$ pointwise. Since $A_{\lambda} \xrightarrow{\mathcal{LE}_{\mathcal{S},d}^-} A$, by Theorem \ref{lowerleftexcesschracterization}, we can find $\lambda_0$ such that $\forall$ $\lambda \geq \lambda_0$ $A_\lambda$ hits $B_d(S'', \epsilon)$ pointwise. Then for any $x \in S$ and $\lambda \geq \lambda_0$, we have \begin{eqnarray*}
	 	d(x, A_{\lambda}) \leq d(x, a_x) + d(a_x,z_x) + d(z_x, A_{\lambda})\\ 
	 	< f(x) + \alpha + \epsilon < g(x).	 	 
	 \end{eqnarray*} So eventually $(A_{\lambda})$ hits $B_d(S,g)$ pointwise. Consequently, by Theorem \ref{Hittypetheoremforlowertsd}, $A_{\lambda} \xrightarrow{\tau_{\mathcal{S},d}^-} A \in CL(X)$.  \end{proof}

\bibliographystyle{plain}
\bibliography{reference_file}
\end{document}